\documentclass[12pt]{article}
\usepackage{amsfonts,amsmath,amssymb,latexsym}
\newtheorem{theorem}{Theorem}

\newtheorem{lemma}[theorem]{Lemma}

\newtheorem{fact}[theorem]{Fact}
\def\bea{\begin{eqnarray}}
\def\eea{\end{eqnarray}}
\def\be{\begin{equation}}
\def\ee{\end{equation}}

\def\Rr{\mathbb{R}}
\def\dd{\mathrm{d}}

\def\<{\langle}
\def\>{\rangle}

\def\proof{\noindent{\it Proof: }}
\newcommand{\po}{{\hspace*{-1ex}}{\bf .  }}

\begin{document}

\title{Killing graphs with prescribed mean curvature and 
Riemannian submersions}
\author{M. Dajczer\thanks{Partially supported by Procad, CNPq and Faperj.} 
\,\,and J. H. de Lira\thanks {Partially 
supported by CNPq and FUNCAP.}}

\date{}
\maketitle

\begin{abstract}
It is proved the existence and uniqueness of  graphs with prescribed
mean curvature in Riemannian submersions fibered by flow lines of a
vertical Killing vector field.
\end{abstract}
\vspace{0.3cm}

{\small \noindent {\bf Keywords:}  Killing graphs, Prescribed 
mean curvature.

\noindent {\bf MSC 2000:} 53C42, 53A10.}

\section{Introduction}

Recent papers devoted to the study of CMC surfaces in certain 
homogeneous three-manifolds are based in the description of these ambient 
spaces as Riemannian submersions over constant curvature model surfaces. 
For instance,
this is the case of \cite{A-R}, \cite{ADR}, \cite{Daniel-1} and
\cite{Daniel-2}. In particular, in \cite{ADR} the authors obtained 
CMC graphs in the
Heisenberg space regarding it as a submersion over $\Rr^2$ 
fibered
by geodesic flow lines of a Killing vector field. The goal in these works is to extend classical results
about CMC surfaces in Euclidean space as well as more recent 
results in nonflat space forms to a more general setting. 

One of the main issues in developing a theory for CMC hypersurfaces 
in general Riemannian ambients is the existence of examples. Methods
which rely mainly on geometric constructions  may fail if the 
ambient space lacks appropriate symmetries or structures. 
However, the problem may be solvable once it is reformulated in analytical terms  as the
existence of CMC graphs for a suitable notion of graph. 
This is the case of Riemannian manifolds
carrying a Killing vector field where the natural notion of
Killing graph has been defined under additional assumptions.

The Dirichlet problem for prescribed mean curvature Killing graphs
in ambient spaces endowed with a Killing field with integrable 
orthogonal distribution was first solved for CMC surfaces in \cite{DR}. 
Then, it was extended in \cite{DHL} to hypersurfaces with prescribed mean curvature function. Under the integrability assumption, the ambient manifold has a warped product structure with one of the factors giving rise to a totally geodesic hypersurface foliation.

In this paper, we consider a generalization of \cite{ADR} 
and \cite{DHL}
to  Riemannian submersion $\pi\colon\,\bar M^{n+1}\to
M^n$ whose vertical fibers are given by flow lines of a 
Killing field. Thus, the normal 
distribution to the Killing field may fail to be integrable. 
Our aim is to show that a natural setting of the Dirichlet 
problem for 
Killing graphs (defined in Section 2) with prescribed mean 
curvature function
in this context is to
consider these as leaves transversal to a solid cylinder of the 
flow lines  that project on a compact domain on the base of the 
submersion. 
Using this approach, we give a unified proof of known and 
completely new existence results in a wide range of ambient 
Riemannian manifolds. Among the ambients for which this paper 
applies, we should mention higher-dimensional 
Heisenberg spaces and odd-dimensional spheres  submersed in 
complex projective spaces.

The existence part of our result is proved using the continuity
 method for
quasilinear elliptic PDE. In order to obtain apriori estimates
essential to this method we use Killing cylinders as barriers. 
Given  a domain $\Omega$ in $M$ with compact closure and boundary
$\Gamma$, the Killing cylinders over $\Gamma$ and
$\bar\Omega$ are, respectively, the subsets $K=\pi^{-1}(\Gamma)$ and
$M_0=\pi^{-1}(\bar\Omega)$. We denote by $H_{\textrm{cyl}}$ the 
inward mean
curvature of $K$ and by $\textrm{Ric}_{\bar M}$ the Ricci 
tensor of $\bar M$.
\medskip

With the above notations we have the following  result.

\begin{theorem}\label{main}\po
Let $\Omega\subset M$ be a domain with compact closure and
$C^{2,\alpha}$ boundary. Suppose that $H_{{\rm cyl}}>0$ and
$
\inf_{\bar M}\textrm{Ric}_{\bar M} \ge -n \inf_\Gamma H_{{\rm cyl}}^2.
$
Let $H\in C^\alpha(\bar\Omega)$ and $\phi\in C^{2,\alpha}(\Gamma)$
be given functions and $\iota\colon\,\bar\Omega \to M_0\subset\bar M$ be a $C^{2,\alpha}$
immersion transversal to the vertical
fibers such that $\pi\circ
\iota=id|_{\bar\Omega}$. If
$$
\sup_\Omega|H|\le \inf_\Gamma H_{\rm{cyl}},
$$
then there exists a unique function $u\in C^{2,\alpha}(\bar\Omega)$
satisfying $u|_\Gamma=\phi$ whose Killing graph $\Sigma$ has mean
curvature $H$.
\end{theorem}

The hypothesis on the existence of an immersion $\iota$ is used
simultaneously to introduce a set of coordinates well suited to the
problem and  to define properly the notion of Killing graph. In terms of
these coordinates, it may be rendered evident that the ambient
metric is stationary.  Moreover, $\iota(\bar\Omega)$ is used as
barrier to producing an initial minimal graph by the direct method
in Calculus of Variations.  In higher-dimensional
Heisenberg spaces there exists a minimal leaf transverse to the flow
lines of the vertical vector field. Thus, in this particular case
there is no need of the hypothesis. By contrast, if we consider the
example of odd-dimensional spheres submersed in the complex
projective spaces, it is not guaranteed that always exist such
minimal graphs with respect to the Hopf fibers. 

We remark
that submersions with totally geodesic fibers constitute an
important example where we may construct initial Killing graphs. In
fact, if we also assume that the Killing cylinder $M_0$ over
$\bar\Omega$ is geodesically complete, then geodesic cones with 
boundary in $K$
and vertex at the mean convex  side of $K$ may be taken, after
smoothing around the vertex,  as initial Killing  graphs. Thus, we
may rule out the hypothesis in this case.

This paper is organized as follows. In Section 2, we fix
notation and made precise the notion of Killing graph. We deduce 
the mean curvature equation and define adapted and basic reference 
frames crucial in the subsequent analysis. In Section 3, 
we present some basic geometry of Killing cylinders. In Sections 4
and 5 we construct  analytical barriers to obtain 
height and boundary gradient estimates.  Section 6 is
devoted to the proof of interior gradient estimates based in 
the technique of normal perturbation of the graph due to  Korevaar \cite{Korevaar}. 
The continuity method and the existence of the minimal initial solution 
are presented in the final section.

\section{Killing graphs}

Let $\pi\colon\,\bar M^{n+1}\to M^n$ be a Riemannian 
submersion such that the leaves of the vertical foliation 
are the trajectories of a nonsingular Killing vector 
field  denoted by $Y\in\mathfrak{X}(\bar M)$. Let 
$\Omega\subset M$ be a $C^{2,\alpha}$ domain  with compact closure.
We assume that the integral curves of $Y$ in
$$
M_0:=\pi^{-1}(\bar\Omega)
$$
are complete lines. Since the hypersurfaces  we work with 
are graphs over $\bar\Omega$ along the integral curves,  
when these curves 
are circles we may pass  to the universal cover of $M_0$ 
without loss of generality.

Let $\iota\colon\,\bar\Omega\to\bar M$ be
an immersion satisfying $\pi\circ\iota=id_{\bar\Omega}$ such that the
hypersurface $\Sigma_0 = \iota(\bar\Omega)$ is transversal to the flow
lines. The initial values for the flow
$\Psi\colon\,\mathbb{R}\times \Sigma_0 \to M_0$ of $Y$ are taken at
$\Sigma_0$, i.e., $\Sigma_0$ corresponds to the level 
hypersurface $s=0$ for the flow parameter $s$. Set 
$\Psi_s =\Psi(s,\,\cdot\,)$. Then, the level 
hypersurfaces $\Sigma_s=\Psi_s (\Sigma_0)$ constitute a foliation of $M_0$ by 
isometric hypersurfaces.

Fix a local reference frame ${\sf v}_1,\ldots, {\sf v}_n$ on
$\bar\Omega$ and set
$$
\sigma_{ij}=\< {\sf v}_i,{\sf v}_j\>.
$$
Let $\bar{\sf v}_1,\ldots,\bar{\sf v}_n$ 
be the corresponding local frame
on $\Sigma_0$, i.e., $\bar{\sf v}_i(p)= \iota_*{\sf v}_i(x)$ if 
$x\in\bar\Omega$ and $p=\iota(x)$. By means of
the flux $\Psi$ we define a local frame at $q=\Psi_s(p)$ in 
$\bar M$ by
$$
\partial_s(q) =\frac{\dd}{\dd s}\, \Psi(s,p) = Y(\Psi(s,p))
=\Psi_*(s,p)\partial_s(p)
$$
and
$$
\bar{\sf v}_i(q)= \big(\Psi_s\circ\iota\big)_*
{\sf v}_i(x).
$$

Let $D_1,\ldots, D_n$ in $\bar M$ denote the basic vector 
fields $\pi$-related to  ${\sf v}_1,\ldots,{\sf v}_n$.
 If $q=\Psi (s,p)$ for $p\in\Sigma_0$, then
$\pi(q)=\pi\circ\Psi (s,p)=\pi(p)$. Therefore,
$$
D_i(q)=\Psi_*(s,p)D_i(p)
$$
since $\Psi_*(s,p) D_i(q)$ is horizontal and
$$
\pi_*(q)\Psi_*(s,p)D_i(p)
=(\pi\circ \Psi)_*(s,p)D_i(p)=\pi_*(p)D_i(p).
$$
That $\pi$ is a Riemannian submersion yields
$$
\< D_i,D_j\>=\< {\sf v}_i,{\sf v}_j\>= \sigma_{ij}.
$$
Setting
$$
D_0 := f^{1/2}\, \partial_s,
$$
we complete a local reference frame
$D_0,D_1,\ldots, D_n$ on $\bar M$ where
$f:=1/|Y|^2$
does not depend on $s$ since $Y$ is a Killing field..

We extend the frame $\bar{\sf v}_1,\ldots,\,\bar{\sf v}_n$ adapted
to the leaves $\Sigma_s$ to a frame $\bar\nabla s,\,\bar{\sf
v}_1,\ldots,\,\bar{\sf v}_n$ in $\bar M$ by  adding the gradient
vector field $\bar\nabla s$ of the function~$s$. Using 
$$
\pi_*(q)\bar{\sf v}_i
=\pi_*(p)\iota_*{\sf v}_i(x)={\sf v}_i(x)=\pi_*(q)D_i
$$
and
$$
1=\partial_s s = \< \bar\nabla s, \partial_s\> =
f^{-1/2}\<\bar\nabla s,D_0\>,
$$
we have that the two frames considered on $\bar M$ are related by
$$
\left\{
\begin{array}{l}
\bar\nabla s = f^{1/2}\, D_0 +\sigma^{ji}D_j(s)\, D_i
\vspace{1ex}\\
\bar{\sf v}_i = \delta_i\, D_0 +  D_i.
\end{array}\right.
$$
The functions  $\delta_i$ are independent of $s$ since
$$
\delta_i = \< \bar{\sf v}_i (q), D_0(q)\> =\<\Psi_{s*}(p)\bar{\sf
v}_i(p),\Psi_{s*}(p)D_0(p)\> =\< \bar{\sf v}_i(p),D_0(p)\>.
$$
Thus,  from
$$
0=\bar{\sf v}_j(s) = \<\bar\nabla s, \bar{\sf v}_j\> =
f^{1/2}\delta_j + D_j(s)
$$
we conclude that  the functions $D_j(s)$ are also independent of $s$.
\vspace{1ex}

The {\it Killing graph\/} $\Sigma=\Sigma_u$ of a function 
$u\in C^2(\bar\Omega)$ is the hypersurface
$$
\Sigma_u=\{\Psi(u(p),p): p\in\Sigma_0\},
$$
where $u$ is seen as a function on $\Sigma_0$ by  taking
$u(p)=u(x)$ when $\pi(p)=x$.
Since $\Sigma$ can also be considered as given by the immersion
$$
\iota_u\colon\,x\in\bar\Omega\mapsto \Psi(u(x),\iota(x)),
$$
its tangent bundle is spanned by the vector fields
\be\label{tan}
(\iota_u)_* {\sf v}_i = {\sf v}_i(u)\,\Psi_s +(\Psi\circ \iota)_*
{\sf v}_i ={\sf v}_i(u)\,\partial_s+\bar{\sf v}_i
=f^{-1/2}{\sf v}_i(u)D_0+\bar{\sf v}_i.
\ee

We may regard $u$ as a function 
in $M_0$ by means of the extension
\be\label{ex}
u(q) = u(x)\;\;\;\mbox{if}\;\;\; \pi(q)=x.
\ee
Thus $D_0(u)=f^{1/2}\partial_s u=0,$
and hence
$$
D_i(u)=\bar{\sf v}_i(u)-\delta_i D_0 (u)
= \bar{\sf v}_i(u)={\sf v}_i(u).
$$
Therefore, we have using (\ref{tan}) that
$$
(\iota_u)_*  {\sf v}_i=(f^{-1/2}D_i(u)+\delta_i)D_0+D_i.
$$
It follows easily that a unit normal vector field to 
$\Sigma$ pointing upwards is 
\be\label{N}
N =\frac{1}{W}(f^{1/2}D_0-\hat u^j D_j),
\ee
where
\be\label{6}
\hat u^j:=\sigma^{ij}D_i(u-s)
\ee
and
$$
W^2 := f+ \sigma_{ij}\hat u^i \hat u^j
=f+ \hat u^i \hat u_i
$$
for $\hat u_i:=\sigma_{ij}\hat u^j$.
Notice that $\hat u^j$ and $W$ can also be seen as 
functions on $M$ since they are independent of $s$.

\subsection{The mean curvature equation}

To compute the mean curvature of $\Sigma$  assume for 
simplicity that the tangent frame ${\sf v}_1,\ldots,{\sf v}_n$ is
orthonormal at $x\in \Omega$.
Hence, the basic frame $D_0,D_1,\ldots,D_n$ is orthonormal 
at points of $\pi^{-1}(x)$.
Thus, 
$$
\< \bar\nabla_{D_0}N,D_0\> =\frac{f^{1/2}}{W}\<
\bar\nabla_{D_0}D_0,D_0\> - \frac{\hat u^j}{W}\<
\bar\nabla_{D_0}D_j,D_0\>
= \frac{1}{W}\< \bar\nabla_{D_0}D_0,\hat u^jD_j\>,
$$
where $\bar\nabla$ denotes the Riemannian connection 
on $\bar M$. We consider on $M$ the vector field
$$
Du:=\hat u^j{\sf v}_j=\sigma^{ij}D_i(u-s){\sf v}_j.
$$
Since $\bar\nabla_{D_0}D_0$ is a horizontal vector field
and $\pi_* (\hat u^jD_j)=\hat u^j{\sf v}_j$, we obtain
$$
\< \bar\nabla_{D_0}N,D_0\>=
\frac{1}{W}\< \pi_*\bar\nabla_{D_0}D_0, Du\>.
$$
By the well-known O'Neill submersion formula \cite{O1}, we have
\be\label{oneill}
\bar\nabla_{D_k}D_j =
\big(\bar\nabla_{D_k}D_j\big)^h+\frac{1}{2}[D_k,D_j]^v.
\ee
Thus, we obtain for $k\ge 1$ that
$$
\< \bar\nabla_{D_k}N,D_k\> 
=-\frac{f^{1/2}}{W}\< \bar\nabla_{D_k}D_k,D_0\>- \<
\bar\nabla_{D_k}\Big(\frac{\hat u^j}{W}D_j\Big),D_k\>
=-\<\nabla_{{\sf v}_k}\frac{Du}{W}, {\sf v}_k\>
$$
where $\nabla$ denotes the Riemannian connection on $M$.
We conclude that
\be\label{div}
nH=-\textrm{div}_{\bar M}N=\textrm{div}_M\,\frac{Du}{W}
-\frac{1}{W}\< \pi_*\bar\nabla_{D_0}D_0, Du\>.
\ee

Denote the covariant derivative in $M$ of 
$Du=\hat u^j{\sf v}_j$ by
$$
\nabla_{{\sf v}_k} Du := \hat u^j_{\textrm{ };k}\,{\sf v}_j
$$
and set
$\hat u_{k;i}:=\sigma_{jk}\hat u^j_{\textrm{ };i}.$
Computing at any point the divergence in (\ref{div})  gives
\begin{eqnarray*}
\textrm{div}_M\,\frac{Du}{W}
\!\!&=&\!\!
\sigma^{ik}\<\nabla_{{\sf v}_i}\frac{Du}{W},{\sf v}_k\>\\
\!\!&=&\!\!\frac{\sigma^{ik}}{W^2}
(W\<\nabla_{{\sf v}_i}Du,{\sf v}_k\>
- {\sf v}_i(W) \<\hat u^j\,{\sf v}_j,{\sf v}_k\>)\\
\!\!&=&\!\!\frac{\sigma^{ik}}{W^3}(W^2\<
\hat u^j_{\textrm{ };i}{\sf v}_j,{\sf v}_k\> -
\frac{1}{2}{\sf v}_i(f + \hat u^l\hat
u_l)   \hat u^j\sigma_{jk})\\
\!\!&=&\!\!\frac{\sigma^{ik}}{W^3}(W^2\sigma_{jk}\hat
u^j_{\textrm{ };i} - \frac{1}{2}({\sf v}_i(f)+2\hat u^l
\hat u_{l;i})\hat u^j\sigma_{jk})\\
\!\!&=&\!\!\frac{1}{W^3}(W^2\sigma^{ik}-\hat u^i\hat
u^k)\, \hat u_{k;i}-\frac{1}{2W^3}\,{\sf v}_i(f) \,\hat u^i.
\end{eqnarray*}
 On the other hand,
\begin{eqnarray*}
{\sf v}_i(f) \!\!&=&\!\! -f^2 {\sf v}_i \< Y,Y\>
=-f^2(\delta_i D_0 \< Y,Y\> +D_i \< Y,Y\>)\\
\!\!&=&\!\!-2f^2\< \bar\nabla_{D_i}Y,Y\>  
= 2f^2\< \bar\nabla_{Y}Y,D_i\>=2f \< \bar\nabla_{D_0}D_0,D_i\>\\
\!\!&=&\!\! 2f \< \pi_*\bar\nabla_{D_0}D_0,{\sf v}_i\>.
\end{eqnarray*}
Thus, the mean curvature equation (\ref{div}) becomes
\be\label{pde}
A^{ik}\hat u_{k;i}
-(f+W^2)\< \pi_*\bar\nabla_{D_0}D_0, Du\>=B,
\ee
where
$$
A^{ik}:=W^2\sigma^{ik}-\hat u^i\hat u^k
\;\;\;\;\mbox{and}\;\;\;\;    B:=nH\, W^3.
$$

We define the operator
$$
\mathcal{Q}[u] =\frac{1}{W^3}(A^{ij} \hat u_{j;i}
-(f+W^2)\<\pi_*\bar\nabla_{D_0}D_0, Du\>).
$$
Therefore, we have shown that $\Sigma$ is a hypersurface with
prescribed mean curvature function $H(x)$ and boundary condition 
$\phi$ if $u$ is a
solution to the~Dirichlet problem

\be\label{dir}
\left\{ \begin{array}{ll}
 \mathcal{Q}[u]=nH\vspace{1ex}\\
u|_\Gamma =\phi
\end{array}\right.
\ee
where $\Gamma=\partial \Omega$.
The boundary of $\Sigma$ is the Killing graph over
$\Gamma$ of $\phi$.

\subsection{A commutation formula}

In this subsection, we give a  commutation formula for 
second covariant derivatives that allow to conclude the 
ellipticity 
of the quasilinear operator defined in the proceeding one.
\vspace{1ex}

Since $Du=\pi_*\bar\nabla(u-s)$ from  (\ref{6})
we obtain using (\ref{oneill}) that
\begin{eqnarray*}
\nabla_{{\sf v}_k} Du \!\!&
= &\!\! \nabla_{{\sf v}_k}\pi_*\bar\nabla(u-s)\\
\!\!&=&\!\!\pi_*(\bar\nabla_{D_k}
\bar\nabla(u-s)-\bar\nabla_{D_k}\<\bar\nabla(u-s),D_0\> D_0)\\
\!\!& = &\!\! \sigma^{il}\pi_* (\<\bar\nabla_{D_k}\,
\bar\nabla(u-s),D_l\>\,D_i-\<\bar\nabla(u-s),D_0\>
\<\bar\nabla_{D_k} D_0,D_l\> \,D_i )\\
\!\!& = &\!\! \sigma^{il} (\bar\nabla^2_{D_k,D_l}\,(u-s)
+(D_0(u)-D_0(s))
\<\bar\nabla_{D_k}D_l, D_0\>)\, \pi_* D_i\\
\!\!& = &\!\!
\sigma^{il}\big(\bar\nabla^2_{D_k,D_l}\,(u-s)
-\frac{1}{2}f^{1/2} \<[D_k,D_l], D_0\> \big)\,{\sf v}_i.
\end{eqnarray*}
Therefore,
\bea\label{dos}
\hat u_{j;k} \!\!& = &\!\!  \sigma_{ji}\hat u^i_{\textrm{
};k}=\sigma_{ji}\sigma^{il}\big(\bar\nabla^2_{D_k,D_l}\,(u-s)
-\frac{1}{2}f^{1/2}
\< [D_k,D_l], D_0\> \big)\nonumber\\
\!\!& = &\!\! \bar\nabla^2_{D_k,D_j}\,(u-s)
+\frac{1}{2}\gamma_{jk}
\eea
where $\gamma_{kj}:=f^{1/2}\< [D_k,D_j],D_0\>$ is skew-symmetric. 
Since the  Hessian is symmetric we conclude that
\be\label{one}
\hat u_{j;k}-\hat u_{k;j}=\gamma_{jk}.
\ee

 Under the convention for $u$ established in (\ref{ex}) we use the 
 standard notation
$$
u_i = D_i(u), \quad u^i = \sigma^{ij}u_j\;\;\;\mbox{and}\;
\;\; u_{i;j}=\<\bar\nabla_{D_i} \bar\nabla u, D_j\>.
$$
Then the matrix $(\hat u_{i;j})$ is
related with the Hessian matrices $(u_{i;j})$ and $(s_{i;j})$~by
\be\label{second}
\hat u_{i;j}= \<\bar\nabla_{D_j}\bar\nabla u,
D_i\>-\<\bar\nabla_{D_j}\bar\nabla s, D_i\>
+\frac{1}{2}\gamma_{ij}
= u_{i;j}-s_{i;j}+\frac{1}{2}\,\gamma_{ij}.
\ee
Hence, the principal part of the mean curvature
equation (\ref{dir}) is  given by
the matrix $(A^{ij})$. This matrix is positive-definite. Indeed, 
we have that
\be\label{eigenvalues}
f|\xi|^2 \le A^{ij}\xi_i\xi_j\le W^2|\xi|^2.
\ee

\section{Killing cylinders}

 The {\it Killing cylinder\/} over $\Gamma$ is  the hypersurface $K=\pi^{-1}(\Gamma)$.  Thus 
$$
K=\{\Psi(s,\iota(x)): s\in\mathbb{R}, \,x\in \Gamma\}
$$
is ruled by the flow lines of $Y$ through 
$\iota(\Gamma)\subset \Sigma_0$.

We denote by $\bar \eta$ the  inward pointing unit vector 
field  normal to  $K$. Clearly, $\bar\eta$ is a basic vector 
field and
$\pi_* \,\bar\eta =\eta$ is the unit normal vector field 
to $\Gamma$ in $M$ pointing inward.
We work with a tangent frame satisfying that
${\sf v}_1=\eta$ and ${\sf v}_2,\ldots, {\sf v}_n$ are 
orthogonal to ${\sf v}_1$. In particular, their horizontal 
lifts $D_i$ verify along $K$ that $D_1=\bar\eta$ and $D_j$, 
$2\le j\le n$, is tangent to $K$. Set
$$
f\<\bar\nabla_YY,\bar\eta\> = \<
\bar\nabla_{D_0}D_0, \bar\eta\>=\kappa,
$$
where $\kappa$ can be seen as a function on $\bar\Omega$.
In fact,
$$
f\<\bar\nabla_YY,\bar\eta\>
=-f\<\bar\nabla_{\bar\eta} Y,Y\>
=\frac{1}{2f}\bar\eta(f).
$$
Thus,
$$
Y(\kappa)=\frac{1}{2f}Y(\bar\eta(f))
=\frac{1}{2f}[Y,\bar\eta](f)=0
$$
since $[Y,\bar\eta]^h=0$ because 
$\pi_*[Y,\bar\eta]^h=[\pi_*Y,\eta]=[0,\eta]$.
Hence, 
$$
n\, H_{\textrm{cyl}} = \sum_{i,j}\sigma^{ij}\<\bar\nabla_{D_i}D_j,
\bar \eta\> +\kappa = \sum_{i,j}\sigma^{ij} \<\nabla_{{\sf v}_i}{\sf
v}_j,\eta\> +\kappa = (n-1)H_\Gamma+\kappa
$$
where $H_\Gamma$ is the mean curvature of $\Gamma$ in $M$.

In the sequel, we deduce some useful properties of the distance
function $d=\textrm{dist}(\,\cdot\,,K)$ from $K$.
We denote by $\Gamma_\epsilon$ and $K_\epsilon$ the level sets 
$d=\epsilon$ in
$M$ and $\bar M$, respectively.  Thus, $\Gamma_\epsilon$ and 
$K_\epsilon$ are
equidistant  from $\Gamma$ and $K$, respectively. It is immediate 
that $K_\epsilon$ is a Killing cylinder over $\Gamma_\epsilon$. Since
$\Gamma$ is assumed to be $C^{2,\alpha}$, the function $d$ is also
$C^{2,\alpha}$ at points of $\Psi(\mathbb{R}\times\Omega_\epsilon)$,
where $\Omega_\epsilon\subset\Omega$ is a tubular 
$\epsilon$-neighborhood
of $\Gamma$ in $M$ for small $\epsilon>0$.

Given
$q\in\Psi(\mathbb{R}\times\Omega_\epsilon)$ we write $q =\exp_{p}\,
d\,\eta$ for some $p\in K$. Hence,
$$
\label{two} |D_1|=|\bar\nabla d|=1.
$$
It follows that
$$
\label{three}
0= \frac{1}{2}D_i\< \bar\nabla d,\bar\nabla d\>
=\<\bar\nabla_{D_i}\bar\nabla d, \sigma^{jk}D_j(d)D_k\>
=d^kd_{i;k}.
$$
We also have
$$
\< \bar\nabla_{D_0}\bar\nabla d,D_0\>
=-\<\bar\nabla_{D_0}D_0,\bar\nabla d\>:=-\kappa_\epsilon
$$
and
$$
\<\bar\nabla_{D_1}\bar\nabla d, D_1\>= \frac{1}{2}D_1 |\bar\nabla d|^2=0.
$$
Therefore,
\be\label{mean}
\Delta d|_{d=\epsilon} 
=-\kappa_\epsilon \sigma^{ij}\<\bar\nabla_{D_i}\bar\nabla d,D_j\>
=  -\kappa_\epsilon-\sigma^{ij}\, b_{ij}^\epsilon
= -nH_{\textrm{cyl}}^\epsilon,
\ee
where $b_{ij}^\epsilon$ are the components of the Weingarten 
operator $A_\epsilon$ and $H_{\textrm{cyl}}^\epsilon$  
the mean curvature of $K_\epsilon$.

\begin{fact}\label{Omega_0}\po
{\em All of the above calculations on
the distance function $d$ remain valid if we replace $\Omega_\epsilon$ 
by the larger  subset $\Omega_0$ in $\Omega$ consisting of the points 
which can be joined to $\Gamma$ by a {\it unique} minimizing geodesic.
It was shown in \cite{LN} that in this set 
$d$ has the same regularity as $\Gamma$. }
\end{fact}

In this paper the ambient Ricci tensor in  direction $v$ is defined~by
$$
\textrm{Ric}_{\bar M} (v)=\sum_{i=1}^n\< \bar R(e_i,v)v,e_i\>,
$$
where $\bar R$ is the curvature tensor in $\bar M$ and
$e_1,\ldots,e_n,v$ is an orthonormal basis. We follow \cite{GT} or
\cite{Sp} and use the result in Fact \ref{Omega_0}  for the proof of
the following result.

\begin{lemma}\label{ricatti}\po
Assume that the  Ricci curvature satisfies  
${\rm Ric}_{\bar M}\ge-n\inf_\Gamma H_{{\rm cyl}}^{\,2}$. 
Let $y_0\in\Gamma$ be the closest point to a given point
$x_0\in\Gamma_\epsilon\subset\Omega_0$. If $H_{{\rm cyl}}>0$, 
then, we have
$$
H_{{\rm cyl}}(\epsilon)|_{x_0}\ge H_{{\rm cyl}}|_{y_0}.
$$
\end{lemma}

\proof  At $d=\epsilon$ and since $D_1$ is the unit speed
of a geodesic, on one hand we have that
\begin{eqnarray*}
 -\frac{\dd}{\dd\epsilon}
\< A_\epsilon D_i,D_j\>
\!\!& =  &\!\!  D_1 \<\bar\nabla_{D_i}D_1,D_j\>
= \<\bar\nabla_{D_1}\bar\nabla_{D_i} D_1,D_j\>
+\<\bar\nabla_{D_i}D_1,\bar\nabla_{D_1} D_j\>\\
\!\!\!\!& = &\!\! -\bar R_\epsilon(D_i,D_j)
-\<\nabla_{[D_i,D_1]}D_1,D_j\>
+\<\bar\nabla_{D_i}D_1,\bar\nabla_{D_1} D_j\>\\
\!\!& = &\!\! -\bar R_\epsilon(D_i,D_j)
-\<\nabla_{D_j}D_1,[D_i,D_1]\>
+\<\bar\nabla_{D_i}D_1,\bar\nabla_{D_1} D_j\>
\end{eqnarray*}
where $\bar R_\epsilon
=\< R(\cdot,D_1)D_1,\cdot\>|_{d=\epsilon}$.
On the other hand,
\begin{eqnarray*}
\frac{\dd}{\dd \epsilon}\< A_\epsilon D_i,D_j\>
\!\!& =  &\!\!
\<\bar\nabla_{D_1}A_\epsilon D_i, D_j\>
+\< A_\epsilon D_i,\bar\nabla_{D_1}D_j\> \\
\!\!&= &\!\!\<(\bar\nabla_{D_1} A_\epsilon)D_i,D_j\>
+\< A_\epsilon D_j, \bar\nabla_{D_1}D_i\>
-\< \bar\nabla_{D_i}D_1,\bar\nabla_{D_1}D_j\>\\
\!\!& = & \!\!\<  A'_\epsilon D_i,D_j\>
-\< \bar\nabla_{D_i}D_1,\bar\nabla_{D_1}D_j\>
-\< \bar\nabla_{D_j}D_1, \bar\nabla_{D_1}D_i\>.
\end{eqnarray*}
Adding the above equations we have the
Ricatti equation
$$
A'_\epsilon - A_\epsilon^2 - \bar R_\epsilon=0,
$$
Taking traces, we obtain
$$
n\frac{\dd}{\dd\epsilon}H_{\textrm{cyl}}^\epsilon = D_1(
\textrm{tr}\,A_\epsilon) =\textrm{tr}\bar\nabla_{D_1} A_\epsilon =
\textrm{tr} \big(A_\epsilon^2 + \bar R_{\epsilon}\big) \ge
n(H_{\textrm{cyl}}^\epsilon)^2 + {\rm Ric}_{\bar M}(D_1).
$$
 From our hypothesis on ${\rm Ric}_{\bar M}$ we have that
$z(d)=H_{\textrm{cyl}}(d)-H_{{\rm cyl}}(y_0)$ satisfies
$$
z'(d)\geq H^2_{\textrm{cyl}}(d)-\inf_\Gamma H_{{\rm cyl}}^2
\geq H^2_{\textrm{cyl}}(d)-H^2_{{\rm cyl}}(y_0)
=(H_{\textrm{cyl}}(d)+H_{{\rm cyl}}(y_0))z(d).
$$
Since $H_{{\rm cyl}}>0$, it follows that $z'(d)\ge c\,z(d)$ in 
some interval $d\in [0,d_0>0]$ for
a constant $c>0$. We obtain  that
$H_{\textrm{cyl}}^\epsilon$ does not decrease with increasing~$d$.
This concludes the proof of the lemma.

\section{The $C^0$ estimate}

In this section, we obtain apriori $C^0$ estimates for
solutions of the Dirichlet problem (\ref{dir}).
\vspace{1ex}

We construct barriers for $u$ in (\ref{dir}) on $\Omega_0$ (see Fact \ref{Omega_0}) by
$$
\label{function} \varphi(x)=\sup_\Gamma \phi + h(d(x))
$$
where $d=\textrm{dist}(\,\cdot\,,\Gamma)$ is regarded as the 
distance from $\Gamma$ on  $M$ and the function $h$ will be 
chosen later.
We work with the frame ${\sf v}_1:=\nabla d,{\sf v}_2, 
\ldots,{\sf v}_n$ 
and the corresponding frame $D_0,D_1,\ldots, D_n$. Thus,
$$
D_i (d)=\<\bar\nabla d, D_i\> = \< D_1,D_i\>
=\< {\sf v}_1,{\sf v}_i\>
= \<\nabla d,{\sf v}_i\>={\sf v}_i(d).
$$

 We have,
$$
\varphi_i = h' d_i \;\;\;\mbox{and}\;\;\;
\varphi_{i;j}=h'' d_i d_j + h' d_{i;j}.
$$
We obtain from (\ref{6}) and (\ref{second}) that
$$
\hat\varphi^j
= \sigma^{ij}(\varphi_i-s_i)=\sigma^{ij}(h'd_i-s_i)
=h'd^i-s^i
$$
and
$$
\hat\varphi_{i;j}=\varphi_{i;j}-s_{i;j}
+\frac{1}{2}\,\gamma_{ij}=h'' d_i
d_j + h' d_{i;j}-s_{i;j}+\frac{1}{2}\,\gamma_{ij}.
$$
Since $\gamma_{ij}$ is skew-symmetric, we have from (\ref{two})
and (\ref{three}) that
\begin{eqnarray*}
A^{ij}\hat \varphi_{j;i}
\!\!&=&\!\! W^2\,\hat\varphi^j_{;i}
-\hat \varphi^i\hat \varphi^j
\hat \varphi_{j;i}\\
\!\!&=& \!\!
W^2(h''+h'd^i_{;i}-s^i_{;i})
-(h'd^i-s^i\big)\big(h'd^j-s^j)
(h''d_id_j+h'd_{j;i}-s_{j;i})\\
\!\!& =&\!\! (W^2-h'^2+2h'\<\bar\nabla d,\bar\nabla s\>
-\<\bar\nabla d,\bar\nabla s\>^2)\,h''
+W^2 h'd^i_{;i}+R
\end{eqnarray*}
where
\be\label{R}
R:=-h's^is^jd_{i;j}-W^2 s^i_{\textrm{ };i}
+(h'^2d^id^j-h'(d^is^j+d^js^i)+s^is^j)s_{i;j}.
\ee
Using that
\be\label{W}
W^2 = f+ \hat\varphi^k\hat\varphi_k
=f+h'^2-2h'\< \bar\nabla d,\bar\nabla s\> +|\bar\nabla s|^2,
\ee
we conclude that
\be\label{final}
A^{ij}\hat \varphi_{j;i}
 =(f+|\bar\nabla s|^2-\<\bar\nabla d,\bar\nabla s\>^2)\, h''
 +W^2 h'd^i_{;i}+R
\ee
where $R$ is a polynomial of second degree in $h'$ and its
coefficients are just functions on $M$.

We have from (\ref{mean}) that
$$
d^i_{;i} =\sigma^{ij}d_{i;j}=\sigma^{ij}\<\bar\nabla_{D_i}\bar\nabla d,D_j\>
=-\sigma^{ij}\, b_{ij}^\epsilon=\kappa_\epsilon-nH_{\textrm{cyl}}^\epsilon.
$$
Since $D\varphi=\pi_*(h'\,\bar\nabla d -\bar\nabla s)$ from (\ref{6}), 
we also have
\be\label{threes}
\<\pi_*\bar\nabla_{D_0}D_0, D\varphi\> =h'\kappa_\epsilon
-\<\bar\nabla_{D_0}D_0,  \bar\nabla s\>.
\ee
Thus, we obtain
$$
W^3\mathcal{Q}[\varphi]   = (f+|\bar\nabla s|^2-\<\bar\nabla
d,\bar\nabla s\>^2)h''
-(f\kappa_\epsilon+nW^2H_{\textrm{cyl}}^\epsilon)h'+R^*
$$
where
\begin{eqnarray*}
R^*:= R + (f+W^2)\< \bar\nabla_{D_0}D_0,\bar\nabla s\>.
\end{eqnarray*}

We choose for (\ref{function}) the test function
$$
h=\frac{e^{CA}}{C}\big(1-e^{-Cd}\big)
$$
where $A>\textrm{diam}(\bar\Omega)$ and $C> 0$ is a  constant 
to be chosen later. Then,
$$
h'=e^{C(A-d)}\;\;\;\mbox{and}\;\;\; h''=-Ch'.
$$
Hence,
$$
\mathcal{Q}[\varphi] \le -(C+\kappa_\epsilon)\frac{fh'}{W^3}
-\frac{h'}{W}nH_{\textrm{cyl}}^\epsilon+\frac{R^*}{W^3}.
$$
Observe that $f/W^2\le 1$. Moreover, as $C\to \infty$ we have that
$1/W\to 0$ and
$$
\frac{h'}{W}=\frac{h'}{(h'^2
-2h'\< \bar\nabla d,\bar\nabla s\>
+|\bar\nabla s|^2+f)^{1/2}}\to 1.
$$
In particular, we have
$$
\frac{R^*}{W^3}\to 0\;\;\;\;\mbox{as}\;\;\;\;C\to \infty.
$$
Choose $C\gg 0$ such that, in particular,  $C+\kappa_\epsilon>0$. 
Using \mbox{$\sup_\Omega|H|\le \inf_\Gamma
H_{\textrm{cyl}}$} and Lemma \ref{ricatti}, we obtain
$$
\mathcal{Q}[\varphi]<-n|H|\le nH.
$$
Thus, one has at points of $\Omega_0$ that 
$$
\mathcal{Q}[\varphi]<\mathcal{Q}[u]=nH,\;\;\;\;\;\;
\varphi|_\Gamma\ge u|_\Gamma.
$$

We now prove that $\varphi\ge u$ on $\bar\Omega$.  By contradiction,
assume that there exist points for which the continuous function
$u^*:=u-\varphi$ satisfies $u^*>0$. Hence $m:=u^*(y)>0$ at a maximum
point $y\in \bar\Omega$ of $u^*$. Choose a minimizing geodesic
$\gamma$ joining $y$ to $\Gamma$ for which the distance
$d=d(y,\Gamma)$ is attained. Thus, $\gamma(t)=\exp_{y_0}t\eta, \;
0\le t\le d$, starts from a point $y_0\in \Gamma$ with unit speed
$\eta$. Since $\gamma$ is minimizing, we have
$d(\gamma(t),\Gamma)=t$ and the function $\varphi$ restricted to
$\gamma$ is differentiable with $ \varphi'(\gamma(t))=e^{C(A-t)}$.
Since the maximum of $u^*$ restricted to $\gamma$ occurs at $t=d$,
i.e., at the point $y$, one has that
$$
u'(\gamma(d))-\varphi'(\gamma(d))=(u^*)'(\gamma(d))\ge 0.
$$
This implies that
$$
\<\nabla u(y),\gamma'(d)\>\ge \varphi'(\gamma(d))=e^{C(A-d)}>0.
$$
In particular $\nabla u(y)\neq 0$, and hence the level hypersurface
$$
S=\{x\in\Omega\cap B_r(y): u(x)=u(y)\}
$$
is regular for small radius $r$. Along $S$ we have
$$
u^* (x)+\varphi(x)=u^* (y)+\varphi(y)\ge u^* (x)+\varphi(y),
$$
and since $\varphi$ is an increasing function of $d(\,\cdot\,,
\Gamma)$ we have  $d(x,\Gamma)\ge d(y,\Gamma)=d$. From this
we conclude that the points in $S$ are at a distance at least $d$
from $\Gamma$.  Since $S$ is $C^2$ it satisfies the interior sphere
condition: there exists a small ball $B_{\varepsilon}(z)$ touching
$S$ at $y$ contained in the side to which  $\nabla u(y)$ and
$\gamma'(d)$ points. Thus, the points of $B_{\varepsilon}(z)$
satisfy $u(x)\ge u(y)$, and hence
$$
\varphi(x)+m\ge u(x)\ge u(y)=\varphi(y)+m,\quad x\in B_\varepsilon
(z),
$$
where in the first inequality we used the definition of $m$. Again
because $\varphi$ is an increasing function of $d$, we have
$d(x,\Gamma)\ge d$ on $B_\varepsilon (z)$ and therefore this ball is
contained in the interior of $\Omega$ far away from $\Gamma$. This
allows us to extend the geodesic $\gamma$ through $B_\varepsilon
(z)$. We claim that the center $z$ of the ball is contained in this
extension. Otherwise, the broken line consisting of $\gamma$ and of
the radius in $B_\varepsilon(z)$ from $z$ to $y$ has length smaller
than {\it a} minimizing geodesic joining $z$ to $y_0\in\Gamma$ (for
a suitable small $\varepsilon$ such a geodesic must cross the level
hypersurface $S$ at a point $x\neq y$ at distance to $\Gamma$
greater than $d$). Thus, if there exists at least two distinct
minimizing geodesics joining $y$ to $\Gamma$, then the point $z$ is
contained in the extension of both geodesics after its intersection
at $y$. Choosing $\varepsilon$ sufficiently small, we see that this
configuration is not possible (the construction we made above
applies to both geodesics). This contradiction implies that the
maximum point $y$ belongs to $\Omega_0$. However, in this case, $u^*
(y)\le 0$, a contradiction. We conclude that $u\le \varphi$
throughout $\bar\Omega$ and therefore $\varphi$ is a continuous
super-solution for the Dirichlet problem (\ref{dir}).

In a similar way, we may construct lower barriers for $u$, that is, 
continuous sub-solutions for (\ref{dir}). It is a standard fact 
that the existence of these barriers implies the desired $C^0$  
apriori estimates. Thus, we have proved the following result.

\begin{lemma}\po\label{height} Under the assumptions of Theorem
\ref{main} there exists a constant $C=C(\Omega,H)$ such that
$$
|u|_{0}\leq C+|\phi|_0
$$
if $u\in C^{2}(\Omega)\cap C^0(\bar \Omega)$ satisfies 
$\mathcal{Q}[u]=nH$ and
$u|_\Gamma=\phi$.
\end{lemma}

\section{Boundary gradient estimates }

In this section our task is to produce apriori gradient estimates
along $\Gamma$ for the Dirichlet problem (\ref{dir}). This is
accomplished by constructing local lower and upper barriers for
$\Sigma$ in a tubular neighborhood of $\Gamma$. \vspace{1ex}

  We construct barriers of the form $w+\phi$ along a tubular neighborhood $\Omega_\epsilon$ of $\Gamma$ as
defined in Section 3. Here, $w=\psi(d(x))$ for some real function
$\psi$ to be chosen and $d=\textrm{dist}(\,\cdot\,,\Gamma)$.
Moreover, the boundary data $\phi$ is extended to a function in 
$\Omega_\epsilon$
along the normal geodesics in a way we make precise later.

We denote
$$
\mathcal{\tilde Q}[u] =\mathcal{Q}[u]-nH.
$$
A simple estimate using 
(\ref{eigenvalues}) and then (\ref{threes}) gives
\bea\label{desi}
\mathcal{\tilde Q}[w+\phi]\!\!&=&\!\!
a^{ij}(x,\nabla w+\nabla \phi)(\hat w_{i;j}+\hat \phi_{i;j})
+b(x,\nabla w +\nabla\phi)-nH\nonumber \\
\!\!&\le&\!\! a^{ij}\hat w_{i;j}
+\frac{1}{W}|\phi|_{2,\alpha} +b-nH,
\eea
where
\be\label{aij}
a^{ij}:=\frac{A^{ij}}{W^3}
=\frac{1}{W}\sigma^{ij}-\frac{1}{W^3}(\hat w^i+\hat\phi^i)(\hat
w^j+\hat \phi^j)\hat w_{i;j}
\ee
and
$$
b = -\frac{f+W^2}{W^3}(\psi'\kappa_\epsilon
+\<\pi_* \bar\nabla_{D_0}D_0, D\phi\>-\<\bar\nabla_{D_0}D_0,\bar\nabla s\>)
$$
since
$\kappa_\epsilon= \<\bar\nabla_{D_0}D_0,\bar\nabla d\>$ and
$$
D(w+\phi)=\sigma^{ij}(\psi'd_j+\phi_j-s_j){\sf v}_j
=D\phi + \pi_*(\psi'\,\bar\nabla d -\bar\nabla s).
$$

 From now on $R_j,\,j\ge 1$, denotes a polynomial
of at most second degree in $\psi'$ whose
coefficients are functions in $M$. As in (\ref{W}) and (\ref{final})
we first obtain,
$$
W^2= f+\psi'^2-2\psi'\<\bar\nabla d,\bar\nabla s
-\bar\nabla \phi\>+|\bar\nabla s-\bar\nabla \phi|^2,
$$
and then 
$$
W^2\hat w^i_{\textrm{ };i}-\hat w^i\hat w^j\hat w_{i;j}
=(f+|\bar\nabla s-\bar\nabla\phi|^2
-\< \bar\nabla d,\bar\nabla s-\bar\nabla\phi\>^2)\psi''
+W^2 \psi'd^i_{;i} +R_1.
$$
 Moreover,
\begin{eqnarray*}
\hat w^i\hat \phi^j\hat w_{i;j}
\!\!&=&\!\!(\psi' d^i-s^i)\hat \phi^j
(\psi''d_i d_j + \psi' d_{i;j}-s_{i;j})\\
\!\!&=&\!\! -\psi'\psi''\<\bar\nabla
d,\bar\nabla s-\bar\nabla \phi\>
+\psi''\<\bar\nabla d,\bar\nabla s\>
\<\bar\nabla d,\bar\nabla s-\bar\nabla \phi\> +R_2
\end{eqnarray*}
and
$$
\hat\phi^i\hat \phi^j \hat w_{i;j}
=\hat\phi^i\hat\phi^j(\psi'' d_i
d_j+\psi'd_{i;j}-s_{i;j}) + R_3 .
$$
Now define
$$
\psi (d) = \mu\ln (1+Kd)
$$
for constants $\mu>0$ and $K>0$ to be chosen later. We have
$$
\psi'  = \frac{\mu K}{1+Kd}\quad\textrm{ and }\quad
\psi''=-\frac{1}{\mu}\psi'^2.
$$
Then using $d^i_{;i}= -n H_{\textrm{cyl}}^\epsilon
+\kappa_\epsilon$  we obtain
$$
W^2\hat w^i_{\textrm{ };i}-\hat w^i\hat w^j\hat w_{i;j}
=-\psi'(n H_{\textrm{cyl}}^\epsilon
-\kappa_\epsilon)W^2 + R_4,
$$
$$
\hat w^i\hat \phi^j\hat w_{i;j}=-\psi'\psi''\<\bar\nabla
d,\bar\nabla s-\bar\nabla \phi\> + R_5
$$
and
$$
\hat\phi^i\hat \phi^j \hat w_{i;j}=R_6.
$$
Since (\ref{aij}) gives
$$
W^3a^{ij} \hat w_{i;j}=W^2\hat w^i_{\textrm{ };i}
-\hat w^i\hat w^j\hat w_{i;j}
-(\hat w^i\hat \phi^j+\hat w^j\hat \phi^i)\hat w_{i;j}
-\hat\phi^i\hat \phi^j\hat w_{i;j},
$$
we now conclude from (\ref{desi}) that
$$
W^3\mathcal{\tilde Q}[w+\phi] \le
-\psi'(nH_{\textrm{cyl}}^\epsilon
-\kappa_\epsilon)W^2
-\frac{2}{\mu}\psi'^3\<\bar\nabla
d,\bar\nabla s-\bar\nabla \phi\>+(b-nH)W^3+R_7.
$$
 From the expressions above for $b$ and $W^2$ it follows that
$$
bW^3 + \psi'\kappa_\epsilon W^2 = R_8.
$$
Hence, we obtain
$$
W^3\mathcal{\tilde Q}[w+\phi] \le
-(n(H+H_{\textrm{cyl}}^\epsilon)
+\frac{2}{\mu}\<\bar\nabla d,\bar\nabla s
-\bar\nabla \phi\>)\psi'^3+R_9.
$$

We choose $\mu$ in such a way that $\mu\to 0$ as $K\to \infty$. Namely,
$$
\mu = \frac{C}{\ln (1+K)}
$$
for some constant $C>0$ to be chosen later. As $K\to\infty$ we have that
$$
\psi'(0)= \frac{CK}{\ln (1+K)}\to +\infty.
$$
It also holds that $\psi'/W\sim 1$ as $K\to \infty$. Thus,
at points of $\Gamma$ the last inequality becomes
$$
W^3\mathcal{\tilde Q}[w+\phi] \le
-(n(H+H_{\textrm{cyl}})
+\frac{2}{\mu}\<\bar\nabla s
-\bar\nabla \phi,\eta\>)\psi'^3+R_9.
$$
We choose the extension of $\phi$ in such a way that at points 
of $\Gamma$ it holds
$$
\<\bar\nabla\phi,\eta\> < \<\bar\nabla s,\eta\>.
$$

Therefore, assuming that $H_{\textrm{cyl}}+H\ge 0$ and choosing $K$ 
large enough, we assure that
$\mathcal{\tilde Q}[w+\phi]<0$ on a
small tubular neighborhood $\Omega_\epsilon$ of $\Gamma$.
Notice that $(w+\phi)|_{\Gamma}=\phi|_{\Gamma}$. Choosing $C$ and $K$ 
large enough we also have that $w+\phi\ge u|_{\Gamma_\epsilon}+\phi$.  
Therefore, $w+\phi$ is a locally defined upper barrier for the Dirichlet 
problem (\ref{dir}). A lower barrier may be constructed in a similar way.
Thus, we have proved the following fact.

\begin{lemma}\po\label{gradbound}
Assume that $u\in C^{2}(\Omega)\cap C^{1}(\bar \Omega)$  satisfies $\mathcal{Q}[u]=nH$ and \mbox{$u|_{\Gamma}=\phi$}. If $|u|$ is bounded in $\bar\Omega$, then 
$$
\sup_{\,\Gamma}|\nabla u|\leq C
$$
by a constant that depends on $|u|_0$.
\end{lemma}

\section{Interior gradient estimates}

\subsection{The prescribed mean curvature case}

In this general case, we adopt ideas from the classical estimate of Korevaar
\cite{Korevaar}. Suppose that the maximum of $|D u|$ is attained at
an interior point, say  $x_0\in \Omega$, where we may assume  that $|D u|\neq 0$ without
loss of generality. Consider a geodesic ball
$B=B(x_0,\rho)\subset\Omega$ centered at $x_0$ with small radius
$\rho\le 1$ so that $|D u|\ge C$ at points of $\bar B$ for some
positive constant $C$. Without loss of generality, we may assume
after a translation along the flow lines of $Y$, if necessary, that
$u<0$ at points of the solid cylinder $\pi^{-1}(\bar B)$.

Let $\eta(x,s)\ge 0$ be a continuous function defined in
$\bar B\times\Rr^-$ that vanishes on $\partial B\times\Rr^-$ and 
is smooth wherever it is positive. Then, let $\bar\Sigma$ be the 
normal geodesic graph over $\Sigma$ defined by
$$
q=\exp_p\epsilon\eta(p) N(p)
$$
where $p\in \Sigma$ is parametrized by $(x,u(x))$. Recall that $N$ given in (\ref{N}) was fixed to be upwards.

For small $\epsilon>0$, we may describe $\bar\Sigma$ as a Killing
graph of some function $\bar u$ defined in $\bar\Omega$. We denote
by $y$ the point in $\Omega$ that maximizes $\bar u -u$. It is clear
that $y\in B$ and that  $D_i\bar u = D_i u$ at this point. From
(\ref{N}) the tangent planes to both graphs have the same slope with
respect to the flow line $\pi^{-1}(y)$ of~$Y$.

We claim that
\be\label{pm}
H_{\bar u}(y)\le H_u(y)
\ee
where $H_u$ and $H_{\bar u}$ denote the mean curvature of $\Sigma$
and $\bar \Sigma$, respectively. In fact, moving $\Sigma$ upward
along the flow lines until the points $(y, u(y))\in\Sigma$ and 
$(y,\bar u(y))\in\bar\Sigma$ coincide, we obtain a tangency point 
for both graphs. Moreover, by the choice of $y$ it is clear that the
translated copy of $\Sigma$ is above $\bar\Sigma$ locally around the
point. Thus, the inequality (\ref{pm}) is consequence of the
comparison principle for the mean curvature PDE.

In analytical terms, the above geometrical reasoning is justified 
in the following way: one has $H_u(x)=\mathcal{Q}[u](x)$ 
and $H_{\bar u}(x)=\mathcal{Q}[\bar u](x)$
since both hypersurfaces are described as Killing graphs.
By construction,  $u=\bar u$ at $\partial B$ and $u\le \bar u$ 
in $\bar B$ (this time, we are not
considering the translation of the geometric proof). If 
$\mathcal{Q}[\bar u]\ge\mathcal{Q}[u]$ in $B$, then the 
analytical comparison principle (cf.\ Thm.\ 10.1 em \cite{GT}) assures that  $\bar u\le u$ in $B$. Thus, this contradiction shows that (\ref{pm})  holds.

It is a well-known fact that since the variation of $\Sigma$ we
consider is along the normal direction, then the mean curvature may be expanded as
\be\label{jacobi}
nH_{\bar u}(\bar x) = nH_u(x)+\epsilon J\eta+O(\epsilon^2),
\ee
where $(x,u(x))$ and $(\bar x, \bar u(\bar x))$ parametrize
correspondent points in $\Sigma$ and $\bar\Sigma$ along the 
same normal geodesic and 
$$
J = \Delta_\Sigma +|A|^2 + \textrm{Ric}_{\bar M}(N,N)
$$
is the Jacobi operator produced by the linearization of the mean 
curvature equation. Here,  $\Delta_\Sigma$ is the Laplace-Beltrami 
operator on $\Sigma$ and $|A|$ denotes the norm of its second 
fundamental form. 

Let $\bar x=y$ for some $x$. It follows from  (\ref{pm}) and  (\ref{jacobi}) that
$$
\epsilon J\eta + O(\epsilon^2) = n(H_{\bar u}(y)
-H_u(x)) \le n(H_u(y)-H_u(x)).
$$
On the other hand, Taylor's expansion of $H_u$ gives
$$
H_u(y) = H_u(x)
+\epsilon\eta H_iT^i + O(\epsilon^2),
$$
where $T^i$ are the components of the horizontal projection of the
normal vector field $N$. Thus, we get at $y$  that
$$
\Delta_\Sigma \eta + |A|^2\eta + \textrm{Ric}(N,N) \eta \le n\eta
H_iT^i + O(\epsilon).
$$
Therefore,
\be\label{delta-1}
\Delta_\Sigma \eta - M\eta\le O(\epsilon)
\ee
for some constant $M>0$ which does not depend on $\eta$.

Next we proceed as in \cite{Korevaar} choosing $\eta=g(\theta(x,s))$
for some real function $g$ to be chosen and a function $\theta$ defined so that
$\Delta_\Sigma \eta$ is large for sufficiently large $|D u(x)|$.
Since $\epsilon$ is chosen small, then (\ref{delta-1})  will give a
contradiction. Observe that $C$ being large implies that the tangent
hyperplanes to $\Sigma$ near $(y,u(y))$ are very sloppy.

That a tangent hyperplane to $\Sigma$ is almost vertical means the
tangential component $\nabla_\Sigma\theta$ of the gradient of
$\theta$ is approximately $\theta_s$. Then, we  define
$$
\theta(x,s)=(Ks+(\rho^2-r^2))^+
$$
for some small constant $K>0$, where $r(x)=\textrm{dist}_M (x_0, x)$ 
is the geodesic distance measured from the center $x_0$ of  $B$ 
and $(\;\cdot\;)^+$ means positive part. We have that 
$0\le\theta\le\rho$.
Since we are assuming height estimates for
$\Sigma$, we may choose $K$ sufficiently small in such a way 
that $\theta>0$ in a neighborhood of $(y,u(y))$ in 
$B\times \Rr^-$. We restrict ourselves to points where $\theta$ is
differentiable. There,
$$
\theta_s =K>0.
$$
Since
\be\label{delta-2}
\Delta_\Sigma \eta =g'' |\nabla_\Sigma \theta|^2
+g'\Delta_\Sigma\theta,
\ee
we have from (\ref{delta-1}) and (\ref{delta-2}) that
\be\label{delta-3}
g''|\nabla_\Sigma \theta|^2 + g' \Delta_\Sigma\theta -Mg\le
O(\epsilon).
\ee
By hypothesis,  the tangent plane of $\Sigma$ at $(y,u(y))$ is not
horizontal. Otherwise, we obtain from (\ref{N}) that $Du(y)=0$.
Let $e$ be the unit vector that gives the steepest ascent 
direction
in the tangent plane of $\Sigma$ at $(y,u(y))$, namely,
$$
e= \frac{1}{W|Du|}(|Du|^2D_0 +f^{1/2}\hat u^jD_j).
$$
Denoting by $\bar\nabla \theta$ the ambient gradient
of $\theta$ and using that $\rho\le 1$, we have
\begin{eqnarray*}
\langle \nabla_\Sigma \theta,e\rangle = \langle\bar\nabla\theta,
e\rangle =\frac{f^{1/2}}{W}\Big(K|Du| 
+\frac{\hat u^jD_j(\theta)}{|Du|}\Big)
\ge\frac{f^{1/2}}{W}(K|Du|-\hat CK  -2),
\end{eqnarray*}
where $\hat C>0$ is a constant independent of $u$ given by the following estimate:
$$
\frac{\hat u^j}{|Du|} D_j (\theta) =\frac{\hat u^j}{|Du|} (KD_j(s)-2r{\sf v}_j(r))
\ge-2 -\hat C K.
$$
Since $K$ and $\hat C$ are independent of $u$ and the
parameter $s$, we may assume that $|Du|> 2/K+\hat C$, and conclude
that
$$
|\nabla_\Sigma \theta|>0.
$$
Finally, for $C_1>0$ large we choose
$$
g(\theta)= e^{C_1 \theta}-1.
$$
It is easily seen that this choice leads to a
contradiction with (\ref{delta-3}). We conclude that $|D u|$ and
therefore $|\nabla u|$ is bounded by some constant which does not
depend on $u$.

\begin{lemma}\po\label{gradglob}
Assume that $u\in C^{3}(\Omega)\cap C^{1}(\bar \Omega)$  
satisfies $\mathcal{Q}[u]=nH$ and \mbox{$u|_{\Gamma}=\phi$}. 
If $u$ is bounded in $\Omega$ and $|\nabla u|$ is
bounded in $\Gamma$, then $|\nabla u|$ is bounded in 
$\Omega$ by a constant that depends
only on $|u|_0$ and $\sup_{\,\Gamma}|\nabla u|$.
\end{lemma}

The usual elliptic regularity results  guarantee that the above estimate  is also true
for a $C^{2,\alpha}$ function (see \cite{GT}).

\subsection{The constant mean curvature case}

In this  case a standard argument works.  In fact, consider the positive function
$$
\Theta:=\< N,Y\> = \frac{f}{\sqrt{f+|Du|^2}}
$$
since a lower estimate of $\Theta$ clearly yields an upper estimate 
for $|\nabla u|$.
 Under the
assumption that $H$ is constant and being $Y$ a Killing field, it is
well-known (cf.\ \cite{ADr}) that $\Theta$ is a Jacobi field, namely,
$J\Theta =0$.
By assumption the Ricci tensor is bounded from below. Thus, since
$\Sigma$ is compact there is a constant $c\geq 0$ such that $ |A|^2
+ \textrm{Ric}_{\bar M}(N,N)\ge -c. $ Thus $\Theta$ is a supersolution 
to the
elliptic operator $\Delta_\Sigma -c$. Hence, the classical minimum
principle states that
$$
\min_\Sigma \Theta \ge \min_{\partial\Sigma}\Theta.
$$
This assures that $|\nabla u|$ is uniformly bounded from above by 
a constant involving the boundary estimates
for $|\nabla u|$.

\section{The proof of the theorem}

In view of the Continuity Method, one must seek for an initial
minimal surface with boundary given by $\Gamma$. This may be
accomplished by defining the sets
$$
\mathcal{C} = \{u\in C^{0,1}(\Omega): u|_\Gamma=0 \}
$$
and, given $k>0$,
$$
\mathcal{C}_k =\{u\in \mathcal{C}: |u|_{0,1}\le k\}.
$$
The hypothesis on the existence of an immersion $\iota\colon\,\bar\Omega\to
\bar M$  assures that the set $\mathcal{C}$ is non-empty
since we may consider the  hypersurface $\iota(\bar\Omega)$ 
as the graph $\Sigma_0$ of the function $u=0$. For the case 
$\kappa=0$, if we assume $M_0$ is geodesically complete, 
the immersion $\iota$ may be obtained as follows:  
we construct a geodesic cone by joining points of a
Killing graph in $K$ over $\Gamma$ to a vertex $p_0$ 
inside $M_0$.
This cone is contained inside $M_0$ since the Killing cylinder
$K$ is mean  convex and $M_0$ is geodesically complete. 
Moreover, it may be smoothed out near the vertex. The
resulting hypersurface may be given as a Killing graph 
since the geodesic cone is always transversal to the 
geodesic vertical fibers.

We then formulate the issue of the existence of a minimal graph
spanning $\iota(\Gamma)$ as the minimization of the functional
$$
\mathcal{I}(u) =\int_\Omega W(x,\nabla u(x))\,\sqrt{\sigma}\dd x ,\quad
u\in\mathcal{C},
$$
where
$$
W= \sqrt{f + \hat u_i \hat u^i}
$$
and the first derivatives of $u$ are taken in a weak sense. Notice
that $f$ and $\hat u_i = u_i -f^{1/2}\delta_i$ do not depend on $u$.
It is clear that $u$ is a critical point of $\mathcal{I}$ if and
only if is a weak solution of the mean curvature equation in
divergence form. Since the principal part of the mean curvature
equation is positive-definite and the coefficients of this equation
(including $H$) do not depend on the function, it follows from
Theorem 11.10 and Theorem 11.11 in \cite{GT} that $\mathcal{I}$  has
a extremum in $\mathcal{C}$. In fact, these theorems require upper
bounds in the Lipschitz norm of the candidates $u\in\mathcal{C}$
which may be obtained from the apriori $C^1$ estimates we derived
earlier.

The $C^{2,\alpha}$ regularity of the minimizer function $u_0$
follows from very general results found in \cite{Morrey}. This
function defines a minimal graph over $\Omega$ with boundary
$\Gamma$.

For the proof of the existence part we apply the well-known continuity method to the family of Dirichlet
problems
\begin{eqnarray*}
\mathcal{Q}_\sigma [u]=n\sigma H,\quad u|_\Gamma = \sigma\phi,
\end{eqnarray*}
where $\sigma\in[0,1]$. The subset $I$ of $[0,1]$ consisting of values of $\sigma$ for which there is a solution is
non-empty since we
have an initial minimal graph spanning the boundary data $\phi$. 
The openness of $I$ is a direct consequence of a standard application of the implicit function theorem since the derivative of $Q_\sigma$  is a linear homeomorphism.
The closedness of $I$ follows from the apriori estimates we had proved and linear elliptic  PDE theory. Thus, the continuity method
assures that $1\in I$.

In order to prove the uniqueness statement, we deduce a kind of flux
formula.  We suppose that there exists a hypersurface $\Sigma'$ in
$M_0$ with $\partial\Sigma'=\Gamma$ and whose mean
curvature is the same as $\Sigma$ at corresponding points in flow
lines. This means that if $x=\pi(p)$ for $p\in \Sigma'$ then the
mean curvature of $\Sigma'$ at $p$ is $H(x)$. Translating $\Sigma'$
we may suppose that $\Sigma$, $\Sigma'$ and a part of the cylinder
$K$ form an oriented cycle which bounds a domain $U$ in $M_0$. Since
$Y$ is tangent to the part of the boundary of $U$ contained in $K$,
we conclude from divergence theorem applied to the field $HY$ in $U$
that
$$
\int_{\Sigma'} \langle HY,N'\rangle=\int_\Sigma \langle HY,N\rangle
$$
where $N$ and $N'$ define respectively the orientations in $\Sigma$
and $\Sigma'$. Applying now the divergence theorem to the
hypersurfaces $\Sigma$ and $\Sigma'$ we obtain that
$$
\int_\Gamma \langle Y,\nu\rangle =\int_{\Gamma'}\langle Y,
\nu'\rangle
$$
where $\nu$ and $\nu'$ are respectively the outward unit co-normals
to $\Gamma$ with respect to $\Sigma$ and $\Sigma'$. This implies
that there exists a point $p$ in $\Gamma$ where $\Sigma$ and
$\Sigma'$ are tangent, that is, where $\nu|_p=\nu'|_p$. Thus, since
$\Sigma'$ is locally a graph near $p$, we conclude from the the
boundary maximum principle that $\Sigma=\Sigma'$.
This concludes the proof of the theorem.

{\renewcommand{\baselinestretch}{1} \hspace*{-20ex}\begin{tabbing}
\indent \= Marcos Dajczer\\
\> IMPA \\
\> Estrada Dona Castorina, 110\\
\> 22460-320 -- Rio de Janeiro -- Brazil\\
\> marcos@impa.br\\
\end{tabbing}}

\vspace*{-2ex}

{\renewcommand{\baselinestretch}{1} \hspace*{-20ex}\begin{tabbing}
\indent \= Jorge Herbert de Lira\\
\> Departamento de Matematica - UFC, \\
\> Bloco 914 -- Campus do Pici\\
\> 60455-760 -- Fortaleza -- Ceara -- Brazil\\
\> jorge.lira@pq.cnpq.br
\end{tabbing}}

\end{document}